\documentclass{amsart}
\catcode`\@=11
\def\@setcopyright{}
\def\serieslogo@{}
\catcode`\@=12

\newcommand{\hrom}[1]{\uppercase\expandafter{\romannumeral#1}}

\def\Q{{\mathbf Q}}
\def\Z{{\mathbf Z}}
\def\C{{\mathbf C}}
\def\R{{\mathbf R}}

\def\Gal{\mathrm{Gal}}
\def\AH{\mathrm{AH}}
\def\Id{\mathrm{id}}
\def\pri{\mathrm{pr}}

\def\End{\mathrm{End}}
\def\Aut{\mathrm{Aut}}
\def\Hom{\mathrm{Hom}}
\def\Lie{\mathrm{Lie}}

\def\GL{\mathrm{GL}}
\def\SL{\mathrm{SL}}
\def\Im{\mathrm{Im}}
\def\dim{\mathrm{dim}}

\def\TORS{\mathrm{TORS}}

\def\Hdg{\mathrm{Hdg}}
\def\det{\mathrm{det}}

\def\k{k}

\def\q{{\mathfrak q}}

\def\Tr{\mathbf{Tr}}

    \def\chil{\bar{\chi}_{\ell}}
\newtheorem{thm}{Theorem}[section]
\newtheorem{lem}[thm]{Lemma}
\newtheorem{cor}[thm]{Corollary}

\theoremstyle{definition}
\newtheorem{defn}[thm]{Definition}

\newtheorem{rem}[thm]{Remark}

\title[Torsion of abelian varieties]
{Torsion of abelian varieties, Weil classes and cyclotomic extensions}

\author[Yuri G. Zarhin]
{Yuri G. Zarhin}

\address{Department of Mathematics, Pennsylvania State University, 
University Park, PA 16802, USA,
\newline
\indent Institute for Mathematical Problems in Biology, 
Russian Academy of Sciences, Push\-chino, Moscow Region, 142292, Russia}
\email{zarhin\char`\@math.psu.edu}
\thanks{Supported by the NSF}

\begin{document}
\maketitle

Let $K\subset\C$ be a field finitely generated over $\Q$, $K(a)\subset \C$ the 
algebraic closure of $K$, $G(K)=\Gal(K(a)/K)$ its Galois group. For each positive integer $m$ we write $K(\mu_m)$ for the subfield of $K(a)$ obtained by adjoining to $K$ all $m$th roots of unity. For each prime $\ell$ we write $K(\ell)$ for the subfield of $K(a)$ obtained by adjoining to $K$ all $\ell-$power roots of unity. We write $K(c)$ for the subfield of $K(a)$ obtained by adjoining to $K$ all roots of unity in $K(a)$. Let $K(ab)\subset K(a)$ be the maximal abelian extension of $K$. The
 field $K(ab)$ contains $K(c)$; if $K=\Q$ then $\Q(ab)=\Q(c)$ (the Kronecker-Weber
 theorem). We write $\chi_{\ell}: G(K)\to \Z_{\ell}^*$ for the cyclotomic
 character defining the Galois action on all  $\ell$-power roots of unity.
 We write $\chil=\chi_{\ell}\bmod{\ell}: G(K)\to \Z_{\ell}^*\to (\Z/\ell\Z)^*$ for
 the cyclotomic character defining  the Galois action on  the $\ell-$th roots of unity.
The character $\chi_{\ell}$ identifies $\Gal(K(\ell)/K)$ with a subgroup
 of  $\Z_{\ell}^*=\Gal(\Q(\ell)/\Q)$. Let $\mu(\Z_{\ell})$ be  the finite cyclic group
$\mu(\Z_{\ell})$ of all roots of unity in $\Z_{\ell}^*$. Its  order is equal
to $\ell-1$ if $\ell$ is odd and $2$ if $\ell=2$.
Let $\Q(\ell)'$ be the subfield of $\mu(\Z_{\ell})-$invariants in
$\Q(\ell)$. Clearly, $\Gal(\Q(\ell)/\Q(\ell)')=\mu(\Z_{\ell})$ and
$\Gal(\Q(\ell)'/\Q)=\Z_{\ell}^*/\mu(\Z_{\ell})$ is isomorphic to $\Z_{\ell}$. 

Let $g$
	be
a positive integer, $X$ a $g-$dimensional abelian variety over $K$. We write
$\End_K(X)$
for the ring of all endomorphisms of $X$ defined over $K$ and $\End^0(X)$
for the finite-dimensional semisimple $\Q-$algebra $\End_K(X)\otimes\Q$.
Its center $F=F_X$ is a field if and only if $X$ is $K-$isogenous to a power of a $K-$simple abelian variety. If so, $F$ is either a totally real number field or a CM-field. We write $\Lie(X)$ for the tangent space to $X$ at the origin. It is the $g-$dimensional $K-$vector space. By functoriality,
$\End^0(X)$ acts faithfully on $\Lie(X)$. We write
$$\Tr_{\Lie(X)}: \End^0(X) \hookrightarrow \End_K(\Lie(X)) \to K\subset \C$$
for the corresponding trace map. The embedding
 $\End^0(X) \hookrightarrow \End_K(\Lie(X))$
 gives rise to  a natural structure of 
(not necessarily faithful) $\End^0(X)\otimes_{\Q}K-$module on $\Lie(X)$.

	The well-known Mordell-Weil-N\'eron-Lang theorem asserts that $X(K)$ is a finitely
generated commutative group. In particular,  its torsion subgroup $\TORS(X(K))$
is finite. Hereafter we will write $\TORS(A)$ for the torsion subgroup of
a commutative group $A$. This implies that $\TORS(X(L))$ is finite for
any finite algebraic extension $L$ of $K$. Mazur \cite{Mazur1} has raised the question of whether the groups $X(K(\ell))$ are finitely generated. In this connection, Serre  (in letters to Mazur, of January 1974) and Imai \cite{Imai} have proved independently that $\TORS(X(K(\ell)))$ is finite for all $\ell$.

Now assume that $L\subset K(a)\subset\C$ is an infinite Galois extension of $K$.
When $L=K(c)$ a theorem of Ribet \cite{Ribet} asserts that $\TORS(X(K(c))$ is finite. The author \cite{ZarhinDuke} has proven that if the center $F$ of $\End^0(X)$ is a direct sum of totally real number fields and $\TORS(X(L))$ is infinite then $L$ contains  infinitely many roots of unity.
On the other hand, Bogomolov (S\'eminaire Delange-Pisot-Poitou, mai 1982, Paris) proved
that $\TORS(X(L))$ is finite if the intersection of $L$ and $K(ab)$ has finite
degree over $K$. For example, if $K=\Q$, we obtain that if $\TORS(X(L))$ is infinite then the intersection of $L$ and $\Q(c)$ has infinite degree over $\Q$. 
The main result of the present paper is the following statement, which deals with essentially non-cyclotomic extensions and may be viewed as a partial improvement of the Bogomolov's result.

\begin{defn}
We say that $X$ and $K$ satisfy hypothesis (H) if they enjoy one of the following
properties:

\begin{enumerate}
\item  There is a discrete valuation $v$  on $K$  such that $X$  has potential purely
multiplicative reduction at $v$;
\item $K$  does not contain a CM-field (e.g., $K\subset \R$);
\item  The Hodge group of $X$  is semisimple.
\item The center $F$ of $\End^0(X)$ is a CM-field and the pair $(X,F)$ is of Weil type, i.e., the $F\otimes_{\Q}K-$module $\Lie(X)$ is free.
\end{enumerate}
\end{defn}

\begin{rem}
It is proven in \cite{SZ} that if an abelian variety
has somewhere a (potential) purely multiplicative reduction then its Hodge
group is semisimple.
\end{rem}

\begin{thm}[Main Theorem]  Let $X$ be a $g$-dimensional abelian variety over
$K$. Assume that $X$  and $K$ satisfy hypothesis (H).
 If the intersection of $L$  and $K(c)$  has finite degree over $K$ then
$\TORS(X(L))$  is finite.
\end{thm}

\begin{rem} 
 If  $L$ is {\sl totally real} then $\TORS(X(L))$ is finite for an arbitrary $X$ \cite{Zhang}.
 We refer to \cite{Ribet}, \cite{ZarhinMatZametki},\cite{ZarhinDuke},
\cite{Wingberg}, \cite{ZarhinMA}, \cite{ZarhinIzv}, \cite{Zhang} for
other results concerning the torsion in infinite extensions.
\end{rem}

The Main Theorem is an immediate corollary of the following  statement.

\begin{thm}
\label{Theorem 1} 
 Let $g$  be a positive integer.  There exists a positive integer $N=N(g)$ depending only on $g$ and enjoying the following properties:

 Let $X$  be a $g-$dimensional abelian variety over $K$  and assume that
 $X$  and $K$ satisfy hypothesis (H).

\sl Then:
\begin{itemize}
\item  Assume that for some prime $\ell$  the $\ell-$primary part of $\TORS(X(L))$  is infinite. Then $K(\ell)$  has finite degree over the intersection  $L \bigcap K(\ell)$  and this degree divides $(N,\ell-1)$ if $\ell$ is odd and divides $2$ if $\ell=2$ respectively. In addition, $L$ contains $\Q(\ell)'$.

\item Let $P=P(X,L)$ be the set of primes $\ell$ such that  $X(L)$ contains a point of order $\ell$.  If $P$  is infinite then for all but finitely many primes $\ell \in P$  the degree $[K(\mu_{\ell}):L\bigcap K(\mu_{\ell})]$  of the field extension $K(\mu_{\ell})/L\bigcap K(\mu_{\ell})$  divides $(N,\ell-1)$.
\end{itemize}
\end{thm}

We will prove Theorem \ref{Theorem 1} in Section 3.

 I would like to thank the MPI f\"ur Mathematik for its hospitality.

\section{Main construction}
 
Let $F$ be the center of $\End_K(X)\otimes\Q$, $R_F=F\bigcap \End_K(X)$ the center of $\End_K(X)$. We put 
$$V_{\Z}=V_{\Z}(X)=H_1(X(\C),\Z), \quad V=V(X)=H_1(X(\C),\Q)= V_{\Z}\otimes\Q.$$
For each nonnegative integer $m$ one may naturally identify the $m$th rational
cohomology group $H^m(X(\C),\Q)$ of $X(\C)$ with $\Hom_{\Q}(\Lambda^m_{\Q}(V(X),\Q)$.
For each prime $\ell$ there are natural identifications
$$X_{\ell}=V_{\Z}/\ell V_{Z}, T_{\ell}(X)=V_{Z}\otimes\Z_{\ell}, V_{\ell}(X)=V(X)\otimes_{\Q}\Q_{\ell}=V_{\Z}\otimes\Q_{\ell}.$$
There is a natural Galois action
$$\rho_{\ell}=\rho_{\ell,X}:G(K)\to \Aut_{\Z_{\ell}}(T_{\ell}(X)) \subset  \Aut_{\Q_{\ell}}(V_{\ell}(X)),$$
induced by the Galois action on the torsion points of $X$ \cite{Serre}.
One may naturally identify the $m$th $\ell-$adic cohomology group $H^m(X_a,\Q_{\ell})$ of $X_a=X\times K(a)$ with
$$\Hom_{\Q_{\ell}}(\Lambda^m_{\Q_{\ell}}(V_{\ell}(X),\Q_{\ell})=\Hom_{\Q}(\Lambda^m_{\Q}(V(X),\Q))\otimes_{\Q}\Q_{\ell}).$$
This identification is an isomorphism of the Galois modules.

Assume now that $F$ is a number field, i.e., $X$ is either simple or is
isogenous over $K$ to a self-product of
	a
simple abelian variety. Let $O_F$
be the ring of integers in $F$. It
is well-known that $R_F$ is a subgroup of finite index in $O_F$. Recall that
for each prime $\ell$ there is a splitting $F\otimes_{\Q}\Q_{\ell}=\oplus F_{\lambda}$
where $\lambda$ runs through the set of prime ideals dividing $\ell$ in $O_F$
and $F_{\lambda}$ is the
completion of $F$ with respect to $\lambda-$adic topology. There is a natural splitting
$V_{\ell}(X)=\oplus V_{\lambda}(X)$
where
$$V_{\lambda}(X)=F_{\lambda} V_{\ell}(X) =V(X)\otimes_F F_{\lambda}.$$
It is well-known that all $V_{\lambda}(X)$ are $G(K)-$invariant $F_{\lambda}-$vector
spaces of dimension
$2\dim(X)/[F:\Q]$. We write $\rho_{\lambda,X}$ for the corresponding
$\lambda-$adic representation 
$$\rho_{\lambda,X}:G(K) \to\Aut_{F_{\lambda}}V_{\lambda}(X)$$
of $G(K)$ \cite{Serre},\cite{RibetA}.
Similarly, for all but finitely many $\ell$
$$R_F/\ell R_F=O_F/\ell O_F = \oplus_{\lambda\mid\ell} O_F/\lambda$$
is a direct sum of finite fields $O_F/\lambda$ of characteristic $\ell$. Also, $X_{\ell}=V_{\Z}/\ell V_{\Z}$ is a free $R_F/\ell R_F=O_F/\ell O_F-$module  of rank $2\dim(X)/[F:\Q]$ and there is a natural splitting
$$X_{\ell}=V_{\Z}/\ell V_{\Z}=\oplus_{\lambda\mid\ell} X_{\lambda}$$
where $X_{\lambda}=(O_F/\lambda) \cdot X_{\ell}.$
Clearly, each $X_{\lambda}$ is a $G(K)-$invariant $O_F/\lambda-$vector space of dimension $2\dim(X)/[F:\Q]$.
We write $\bar{\rho}_{\lambda,X}$ for the corresponding modular
representation
$$\bar{\rho}_{\lambda,X}:G(K) \to\Aut_{O_F/\lambda}X_{\lambda}$$
of $G(K)$ \cite{RibetA}.
Let $d$ be a positive integer and assume that there exists a non-zero
$2d-$linear  form $\psi \in \Hom_{\Q}(\otimes^{2d}_{\Q} V(X),\Q)$, enjoying the following properties. 
\begin{enumerate}
\item  For all $f\in F; v_1, \ldots v_{2d}\in V(X)$
$$\psi(f v_1,v_2,\ldots ,v_{2d})=\psi(v_1,fv_2,\ldots ,v_{2d})=\cdots =
\psi(v_1,v_2,\ldots ,fv_{2d}).$$

\item
For any prime $\ell$ let us extend $\psi$ by $\Q_{\ell}-$linearity to the
non-zero multilinear form
$\psi_{\ell} \in \Hom_{\Q_{\ell}}(\otimes^{2d}_{\Q_{\ell}} V_{\ell}(X),\Q_{\ell})$.
Then for all $\sigma \in G(K); v_1, \ldots v_{2d}\in V_{\ell}(X)$
$$\psi_{\ell} (\sigma(v_1),\sigma(v_2),\ldots,\sigma(v_{2d}))= 
\chi_{\ell}^d(\sigma)\psi_{\ell}(v_1,v_2,\ldots ,v_{2d}).$$
\end{enumerate}
We call such a form {\sl admissible} or $d-${\sl admissible}.

\vskip .5cm

{\bf Example.} Let us assume that $F$ is a {\sl totally real} number field. If
$\mathcal L$ is an invertible sheaf on $X$ defined over $K$ and algebraically
non-equivalent to zero then one may associate to $\mathcal L$ its first Chern
class 
$$c_1({\mathcal L})\in H^2(X(\C),\Q)=\Hom_{\Q}(\Lambda^2_{\Q}(V(X),\Q).$$
The well-known properties of  Rosati involutions and Weil pairings imply
that $c_1({\mathcal L})$ is $1-$admissible (see p.~237 of \cite{MumfordAV}, especially the last sentence and Section 2 of \cite{SZMathZ}).

\vskip .5cm

There exists a unique  $F-2d-$linear form
$\psi_F\in \Hom_F(\otimes^{2d}_F V(X),F)$ such that
$$\psi={\Tr}_{F/\Q}(\psi_F).$$
Multiplying $\psi$ by a sufficiently divisible positive integer, we may and
will assume that the restriction of $\psi_F$ to
$V_{\Z}\times \cdots V_{\Z}$ takes on values in $R_F$. Let
$\Im(\psi_F)$ be the additive subgroup of $R_F$ generated by values of
$\psi_F$ on $V_{\Z}\times \cdots V_{\Z}$ takes on values in $R_F$. Let
$\Im(\psi_F)$ be the additive subgroup of $R_F$ generated by values of
$\psi_F$ on $V_{\Z}\times\cdots V_{\Z}$. Clearly, $\Im(R_F)$ is a
subgroup of finite index in $R_F$ that is an ideal.
Notice that for all but finitely many primes $\ell$ 
$$O_F=R_F/\ell R_F, \Im(\psi_F)=R_F/\ell R_F.$$
Let us extend $\psi_F$ by $F_{\lambda}-$linearity to the
{\sl non-zero} multilinear form
$$\psi_{F,\lambda} \in \Hom_{F_{\ell}}(\otimes^{2d}_{F_{\lambda}} V_{\lambda}(X),F_{\lambda}).$$
Then
$$\psi_{F,\lambda}(\sigma(v_1),\sigma(v_2),\ldots,\sigma(v_{2d}))= 
\chi_{\ell}^d(\sigma)\psi_{F,\lambda}(v_1,v_2,\ldots ,v_{2d})$$ 
for all $\sigma \in G(K); v_1, \ldots v_{2d}\in V_{\lambda}(X)$.

Similarly, for all but finitely many $\ell$ the form $\psi_F$ induces a 
 non-zero multilinear form
$$\psi_{F}^{(\ell)} \in \Hom_{R_F/\ell R_F}(\otimes^{2d}_{R_F/\ell R_F} X_{\ell},R_F/\ell R_F)$$
enjoying the following properties:
\begin{itemize}
\item
The subgroup of $R_F/\ell R_F$ generated by all the values of $\psi_{F}^{(\ell)}$ coincides with $R_F/\ell R_F$;

\item For all $\sigma \in G(K); v_1, \ldots v_{2d}\in X_{\ell}$
$$\psi_{F}^{(\ell)}(\sigma(v_1),\sigma(v_2),\ldots,\sigma(v_{2d}))= 
\chil^d(\sigma)\psi_{F}^{(\ell)}(v_1,v_2,\ldots ,v_{2d}).$$
\end{itemize}

This implies that for all but finitely many $\ell$ the restriction of $\psi_{F}^{\ell}$
to $X_{\lambda}$ defines a {\sl non-zero} multilinear form
$$\psi_{F}^{(\lambda)} \in \Hom_{O_F/\lambda}(\otimes^{2d}_{O_F/\lambda} X_{\ell},O_F/\lambda)$$
enjoying the following property:

$$\psi_{F}^{(\lambda)}(\sigma(v_1),\sigma(v_2),\ldots,\sigma(v_{2d}))= 
\chil^d(\sigma)\psi_{F}^{(\lambda)}(v_1,v_2,\ldots ,v_{2d})$$
for all $\sigma \in G(K); v_1, \ldots v_{2d}\in X_{\lambda}$.

\begin{rem} Using the K\"unneth formula for $X_a^{2d}$, one may view $\psi_{\ell}$ as a
 Tate cohomology class on $X_a^{2d}$. If $\psi$ is skew-symmetric then $\psi_{\ell}$ is a Tate
 cohomology class on $X_a$.
\end{rem}

\begin{thm}  Assume that the center $F$  of $\End^0 X$
 is a field
and there is a $d-$admissible form
$\psi$  on $X$.  Let $\ell$  be a prime and assume that the
$\ell-$ torsion
in $X(L)$  is infinite. If $L^{(\ell)}$  is the intersection of
$L$ and $K(\ell)$
 then the field extension $K(\ell)/L^{(\ell)}$  has finite degree dividing $(d,\ell-1)$ if $\ell$ is odd and dividing $2$ if $\ell=2$ respectively. In addition, $L$ contains $\Q(\ell)'$.

\end{thm}

\begin{proof} As explained in (\cite{ZarhinMA}, 0.8, 0.11) the assumption
that the $\ell-$torsion in $X(L)$ is infinite means that there exists a place
$\lambda$ of F, dividing $\ell$
such that the Galois group $G(L)$ of $L$ acts trivially on $V_{\lambda}(X)$. Since $\psi_{F,\lambda}$ is not identically zero, we conclude that
$$\chi_{\ell}^d(\sigma)=1 \quad \forall \sigma \in G(L) \subset G(K).$$
We write $G'$ for the kernel of $\chi_{\ell}^d$. We have
$G(L)\subset G'\subset G(K)$.

Recall that the kernel of $\chi_{\ell}:G(K) \to \Z_\ell^*$ coincides with the Galois group
 $G(K(\ell))$ of $K(\ell)$ and $\chi_{\ell}$ identifies $\Gal(K(\ell)/K)$ with a subgroup of $\Z_\ell^*=\Gal(\Q(\ell)/\Q)$. Since the torsion subgroup of $\Z_\ell^*$ is the cyclic group $\mu(\Z_{\ell})$ of order $\ell-1$ if $\ell$ is odd and of order $2$ if $\ell=2$, $G'$ coincides with the kernel of $(\chi_{\ell})^{d'}$ with $d'=(d,\ell-1)$ if $\ell$ is odd and $d'=(d,2)$ if $\ell=2$ respectively. This implies that the field $K'=K(a)^{G'}$ of $G'-$invariants is a subfield of $K(\ell)$ and $[K(\ell):K']$ divides $d'$, since $\chi_{\ell}$ establishes an isomorphism between $\Gal(K(\ell)/K')$ and 

$$\{s \in \Im(\chi_{\ell})\subset \Z_{\ell}^*\mid s^{d'}=1\}
\subset \{s \in \mu(\Z_{\ell})\mid s^{d'}=1\}.$$
 Now it is clear that $K'\subset L$, since $G(L) \subset G'=G(K')$.
 It is also clear that $K(\ell)/K'$ is a cyclic extension of degree dividing $d'$.

In order to prove the last assertion of the theorem, notice that
$\Gal(K(\ell)/K) \subset \Gal(\Q(\ell)/\Q)=\Z_{\ell}^*$ and the finite subgroup $\Gal(K(\ell)/K')$ of $\Gal(K(\ell)/K)$ sits in $\mu(\Z_{\ell})\subset\Z_{\ell}^*$. Since $\mu(\Z_{\ell})=\Gal(\Q(\ell)/\Q(\ell)')$, $\Q(\ell)'\subset K'$. Since $K'\subset L$, $\Q(\ell)'\subset L$.
\end{proof}

\begin{thm}  Assume that the center $F$  of $\End^0 X$
 is a field and there is a $d-$ admissible form $\psi$  on
$X$.  Let $S$  be an infinite set of primes $\ell$   such that
for all but finitely many $\ell\in S$  the $\ell-$torsion in
$X(L)$  is not
zero. Then for all but finitely many $\ell\in S$  the field extension 
$K(\mu_{\ell})/K(\mu_{\ell})\bigcap L$  has degree dividing $(d,\ell-1)$.
\end{thm}

\begin{proof} 
For all but finitely
many $\ell$ the $G(K)-$module
$X_{\ell}$ is semisimple and the centralizer of $G(K)$ in $\End(X_{\ell})$
coincides with
$\End_K(X)\otimes \Z/\ell\Z$. This assertion was proven in
\cite{ZarhinInv} for number fields $K$; the proof is based on results of Faltings
\cite{Faltings1}. (See \cite{MW} for an
effective version.) However, the same proof works for arbitrary finitely
generated fields $K$, if one uses results of \cite{Faltings2},
generalizing the results of \cite{Faltings1}. 
Clearly,  for all but finitely many $\ell$ the center of $\End_K(X)\otimes \Z/\ell\Z$ 
coincides 
with $R_F/\ell R_F=O_F/\ell O_F$. Applying Theorem 5f of \cite{ZarhinDuke}
to $G=G(K), G'=G(L), H=X_{\ell}, D=\End_K(X)\otimes \Z/\ell\Z,
R=F_F/\ell R_F$, we conclude that for all but finitely many $\ell \in S$ there exists
 $\lambda\mid\ell$ such that $G(L)$ acts trivially on $X_{\lambda}$. Using the Galois
  equivariance of the non-zero form $\psi_{F}^{(\lambda)}$, we conclude that
  for all but finitely many $\ell\in S$ the character
  $\chil^d$ kills $G(L)$. We write $G'$ for the kernel of $\chil^d$. We have
$G(L)\subset G'\subset G(K)$.

Recall that the kernel of $\chil:G(K) \to (\Z/\ell \Z)^*$ coincides with
$G(K(\mu_{\ell}))$ and
$(\Z/\ell \Z)^*$ is a cyclic group of order $\ell-1$.  This implies that
the field $K'=K(a)^{G'}$ of
$G'-$invariants is a subfield of $K(\mu_{\ell})$ and $[K(\mu_{\ell}):K']$
divides $(\ell-1,d)$, since
$\chil$ establishes an isomorphism between $\Gal(K(\mu_{\ell})/K')$ and $\{s \in \Im(\chil)\subset (\Z/\ell \Z)^*\mid s^d=1\}$. One has only to notice that $K'\subset L$, since $G(L) \subset G'=G(K')$.

\end{proof}

\begin{cor}  Assume that the torsion subgroup of  $X(L)$
 is infinite. Then the intersection of $L$ and $K(c)$  has infinite degree over $K$.
\end{cor}

\begin{proof} Indeed, either there is a prime $\ell$ such that the $\ell-$torsion in $X(L)$ is infinite or for infinitely many primes $\ell$ the $\ell-torsion$ in $X(L)$ is not zero. Now, one has only to apply the previous two theorems.
\end{proof}

\section{Weil classes and admissible forms}
Suppose $A$ is an abelian variety defined over $K$, $\k$ is a CM-field,
$\iota : \k \hookrightarrow \End_K^0(A)$
is an embedding, and $C$ is an algebraically
closed field containing $K$ (for instance, $C=\C$ or $C=\bar{\Q}$). Let $\Lie(A)$ be the tangent space of $A$ at the origin,
an $K$-vector space. If $\sigma$ is an embedding of $\k$ into $C$, let
$$n_\sigma = \dim_C\{t \in \Lie(A)\otimes_K C : 
\iota(\alpha)t = \sigma(\alpha)t {\text{ for all }} \alpha \in \k\}.$$
Write ${\bar \sigma}$ for the composition of $\sigma$ with the involution
complex conjugation of $\k$. 

Recall that a triple $(A,\k,\iota)$ is {\em of Weil type} if $A$ is an abelian
variety over an algebraically closed field $C$ of 
characteristic zero, $\k$ is a CM-field, and 
$\iota : \k \hookrightarrow \End^0(A)$ is an embedding, such that
 $n_\sigma = n_{\bar \sigma}$ for all embeddings $\sigma$ of $\k$ into $C$.

It is known \cite{SZMathZ} that $(A,\k,\iota)$ is of Weil type if and only
if $\iota$ makes
$\Lie(A) \otimes_K C$ into a free $\k \otimes_\Q C$-module 
(see p.~525 of \cite{Ribet} for the case where $\k$ is an imaginary quadratic field). 
Now, assume that  $A=X$ and the image $\iota(k)$ contains the center $F$ of
 $\End_K(X)\otimes\Q$ (for instance, $F=k$). Notice that in the case of Weil type the
  degree $[k:\Q]$ divides $\dim(A)$. In particular, $\dim(A)$ is even.

Our goal is to construct an admissible form, using a triple $(A,\k,\iota)$  of Weil type. 

Recall that  the degree $[k:\Q]$ divides $\dim(A)$, put $d=\dim(X)/[k:\Q]$
and consider the space of Weil classes (\cite{WeilHodge}, \cite{Deligne}, \cite{SZMathZ})
$$W_{k,X}=\Hom_k(\Lambda_k^{2d} V(X),\Q(d)) \hookrightarrow \Hom_{\Q}(\Lambda_{\Q}^{2d}V(X),\Q(d))=H^{2d}(X(\C),\Q)(d).$$
Clearly, $W_{k,X}$ carries a natural structure of one-dimensional $k-$vector
space. If fix an isomorphism of one-dimensional $\Q-$vector spaces $\Q \cong
\Q(2d)$ then one may naturally
identify $\Hom_{\Q}(\Lambda_{\Q}^{2d}V(X),\Q(d))$ with $\Hom_{\Q}(\Lambda_{\Q}^{2d}V(X),\Q)$
and $W_{k,X}$ can 
be described as the space of all $2d-$linear skew-symmetric form
$\psi \in \Hom_{\Q}(\Lambda^{2d}_{\Q} V,\Q)$ with
$$\psi(f v_1,v_2,\ldots ,v_{2d})=\psi(v_1,fv_2,\ldots ,v_{2d})=\cdots =
\psi(v_1,v_2,\ldots ,fv_{2d})$$
for all $f\in F; v_1, \ldots v_{2d}\in V(X)$.

Since $(X,k.\iota)$ is of Weil type, all elements of $W_k$ are Hodge classes by Proposition 4.4 of \cite{Deligne}. Therefore, by Theorem 2.11 of \cite{Deligne} they  must be also {\sl absolute Hodge cycles}; cf. \cite{Deligne}. 

\begin{lem} Let $\mu_k$ be the finite multiplicative group of all roots of
unity in $k$.  There is a continuous
character $\chi_{X,k}:G(K) \to \mu_k \subset k^*,$
 enjoying the following properties:

 For each prime $\ell$  the subgroup
$$W_k \subset W_{k,X}\otimes_{\Q}\Q_{\ell}\subset H^{2d}(X(\C),\Q)(d)\otimes_{\Q}\Q_{\ell}
=H^{2d}(X_a,\Q_{\ell})(d)$$
 is $G(K)-$stable and the action of $G(K)$  on $W_k$ 
is defined via the character
$$\chi_{X,k}:G(K) \to \mu_k \subset k^* =\Aut_k(W_{k,X}).$$
\end{lem}
\begin{proof} Since all elements of $k$ are endomorphisms of $X$ defined over $K$, 
it follows easily that $W_{k,X}\otimes_{\Q}\Q_{\ell}$ is $G(K)-$stable and $G(K)$ acts 
on $W_{k,X}\otimes_{\Q}\Q_{\ell}$ via a certain character
$\chi_{X,k,\ell}:G(K) \to  [k\otimes_{\Q}\Q_{\ell}]^*=
\Pi_{\lambda\mid \ell}k_{\lambda}^*.$

Let us consider the $\Q-$vector subspace 
	$$C^d_{\AH}(X)\subset H^{2d}(X(\C),\Q)(d)\subset H^{2d}(X_a,\Q_{\ell})(d)$$
 of absolute Hodge classes.  Then $C^d_{\AH}(X)$ is $G(K)-$stable and the action of  $G(K)$ on $C^d_{\AH}(X)$ does not depend on $\ell$ and factors through a finite  quotient; cf. \cite{Deligne}, Prop. 2.9b. Since $W_{k,X}$ consists of Hodge classes  and $X$ is an abelian variety, all Weil classes are absolute Hodge classes, i.e,
	$W_{k,X}\subset C^d_{\AH}(X),$
\cite{Deligne}, Th. 2.11. This implies easily that
the subgroup $\Im(\chi_{X,k,\ell})$ is  finite and contained in $k^*$, since the
intersection of $W_{k,X}\otimes_{\Q}\Q_{\ell}$ and $C^d_{\AH}(X)$ coincides with
$W_{k,X}$. (In fact, $W_{k,X}$ coincides even with the intersection of $W_{k,X}\otimes_{\Q}\Q_{\ell}$  and $H^{2d}(X(\C),\Q)(d)$.) This implies also that $\chi_{X,k,\ell}$ does not depend on the choice of $\ell$. So, we may view $\chi_{X,k,\ell}$ as the continuous homomorphism
$$\chi_{X,k}:=\chi_{X,k,\ell}:G(K) \to \mu_k \subset k^*,$$
which does not depend on the choice of $\ell$.
\end{proof}

Let $r$ be the order of the finite group $\Im(\chi_{X,k})$. Clearly, $r$ divides the order of $\mu_k$. Let us put $Y=X^r$ and consider the K\"unneth chunk
$$H^{2d}(X(\C),\Q)(d)^{\otimes r} \subset H^{2dr}(X(\C)^r,\Q)(dr)=H^{2dr}(Y(\C),\Q)(dr)$$
of the $2dr$th rational cohomology group of $Y$. One may easily check that the tensor power
$$W_{k,X}^{\otimes r}\subset H^{2d}(X(\C),\Q)(d)^{\otimes r} \subset H^{2dr}(X(\C)^r,\Q)(dr)=H^{2dr}(Y(\C),\Q)(dr)$$
coincides with the space $W_{k,Y}$ of Weil classes on $Y$ attached to the ``diagonal" embedding
$$k \to \End^0(X) \subset \End^0(X^r)=\End^0(Y).$$
Since the centers of $\End^0(X)$ and $\End^0(X^r)$ coincide, the image of $k$ in $\End^0(Y)$ does contain the center of $\End^0(Y)$.

One may easily check that $G(K)$ acts on $W_{k,Y}=W_{k,X}^{\otimes r}$ via the character $\chi_{X,k}^r$, which is trivial, i.e., $W_{k,Y}$ consists of $G(K)-$invariants.

Let us fix an isomorphism of one-dimensional $\Q-$vector spaces $\Q \cong \Q(2dr)$ and
choose a {\sl non-zero} element
$$\psi \in W_{k,Y} \subset H^{2dr}(Y,\Q)(dr)=\Hom_{\Q}(\Lambda^{2dr}_{\Q}(V(Y),\Q).$$
Then a skew-symmetric $2dr-$linear form $\psi$ is admissible.

Applying to $\psi$ the theorems of the previous section, we obtain the following
statement, which
implies the case 4 (in the hypothesis (H)) of Theorem \ref{Theorem 1}.

\begin{thm}  Assume that the center $F$  of $\End^0 X$
 is a CM-field and $(X,F, \Id)$  is of Weil type. Let us put
$d= \#(\mu_F) \times 2 \dim(X)/[F:\Q] \in \Z_{+}.$
 Let $L$  be an infinite Galois extension of $K$.

\begin{enumerate}
\item  Let $\ell$ be a prime such that the $\ell-$torsion
in $X(L)$  is infinite. Let $L^{(\ell)}$  be the intersection of
$L$  and $K(\ell)$.
 Then the field extension $K(\ell)/L^{(\ell)}$  has finite degree dividing $(d,\ell-1)$ if $\ell$ is odd and dividing $2$ if $\ell=2$ respectively.  In addition, $L$ contains $\Q(\ell)'$.

\item  Let $S$  be the set of primes $\ell$  such that
 $X(L)$  contains a point of order $\ell$  and assume that $S$
  is infinite. Then for all but finitely many $\ell\in S$  the field
 extension
 $K(\mu_{\ell})/K(\mu_{\ell})\bigcap L$  has degree dividing $(d,\ell-1)$.
\end{enumerate}
\end{thm}

\begin{rem} Since $[F:\Q]$ divides $2\dim(X)=2g$, one may easily find an explicit
 positive integer $M=M(g)$, depending only on $g$ and divisible by 
 $\#(\mu_F) \times 2 \dim(X)/[F:\Q]$
\end{rem}

\section{Proof of Theorem \ref{Theorem 1}}

We may and will assume that $X$ is $K-$simple. Then the center $F$ of $\End^0
X$ is either a totally real number field or a CM-field. If $F$ is totally
real then the assertion of
Theorem \ref{Theorem 1} is proven in \cite{ZarhinDuke} with $N=1$. So, further
we assume that $F$ is a
CM-field. We also know that the assertion of Theorem \ref{Theorem 1} is true
when $(X,F)$ is of Weil type
(Case 4 of Hypothesis (H)).

\subsection{ Cases 1 and 3 of Hypothesis (H)} Enlarging $K$ if necessary, we may and will assume
that $X$ is absolutely simple
and has semistable reduction. Then, the results of \cite{SZ} imply that in
both cases $\Hdg_X$ is semisimple.
This means that $(X,F,\Id)$ is of Weil type (cf. for instance \cite{SZ}).
Now, one has only to apply the
result of the previous section with $d=\#(\mu_F)  \times 2\dim(X)/[F:\Q]$
and get the assertion of Theorem \ref{Theorem 1} with $N=M(g)$.

\subsection{Case 2 of Hypothesis (H)} We know that the assertion of the theorem is true if $(X,F,\Id)$ is  of Weil type. So,  we may assume that $(X,F,\Id)$ is not of Weil type.

	Let us consider the trace map
$$\Tr_{\Lie(X)}: F \subset \End^0(X)\hookrightarrow\End_K(\Lie(X)) \to K\subset \C.$$
Our assumption means that the image $\Tr_{\Lie(X)}(F)$ is not contained in $\R$. On the other hand, let us fix an embedding of $F$ into $\C$ and let $L$ be the normal closure of $F$ into $\C$. Clearly, $L$ is a CM-field, containing $\Tr_{\Lie(X)}(F)$. Since $\Tr_{\Lie(X)}(F)\subset K$, the intersection $L\bigcap K$ contains an element, which is not totally real. Since any subfield of a CM-field is either totally real or CM, the field $L\bigcap K$ is a CM-subfield of $K$.

\begin{rem} If $K$ is a number field  not containing a CM-field, one may give another proof, using theory of abelian $\lambda-$adic representations
\cite{Serre} instead of Weil/Hodge classes. The crucial point is that in this case the Serre's tori $T_{\mathfrak m}$ are isomorphic to the multiplicative group ${\mathbf G}_m$ \cite{Serre}, Sect. 3.4.
\end{rem}

\begin{cor}  Let $X$ be a $K-$simple abelian variety of odd dimension. Assume that $K$  does not contain a CM-subfield (e.g., $K\subset \R$).  If $X(L)$  contains infinitely many points of finite order then $L$  contains infinitely many roots of unity.
\end{cor}

\begin{proof} In the case of the totally real center $F$ this assertion is proven in (\cite{ZarhinDuke}, Th.6, p. 142) without restrictions on the dimension. So, in order to prove Corollary, it suffices to check that $F$ is not a CM-field.

	Assume that $F$ is a CM-field. Since $\dim(X)$ is odd, $(X,F,\Id)$ is not
        of Weil type. Now, the arguments, used in the proof of
        Case 2 imply that $K$ contains a CM-subfield. This gives us a desired contradiction.
\end{proof}

\begin{rem} The assertion of Corollary cannot be extended to the even-dimensional
case. In Section \ref{roots} we give an explicit counterexample.
\end{rem}
\begin{rem} Let $X$ be a $g-$dimensional abelian variety that is not necessarily $K-$simple  and let $F$ be the center of $\End^0(X)$. Assume that
$$\Tr_{\Lie(X)}(F)\subset\R.$$
Then the assertion of Theorem \ref{Theorem 1} holds true for $X$.
Indeed, if $Y$ is a $K-$simple abelian subvariety of $X$ and $F_Y$ is the center of $\End^0(Y)$ then one may easily check that either $F_Y$ is a totally real number field or $(Y,F_Y,\Id)$ is of Weil type.
\end{rem}

\section{Example}
\label{roots}
In this section we construct an abelian surface $X$ over $\Q$ and a Galois
extension $L$ of $\Q$ such that $L$ contains only
 finitely many roots of unity but $X(L)$ contains infinitely
many points of finite order.
Of course, the intersection of $L$ and $\Q(c)$  is of infinite degree over $\Q$.

\subsection{} Let $E$ be an elliptic curve over $\Q$ without complex multiplication
(e. g., $j(E)$ is not an
integer). Let us put
$$Y=\{(e_1,e_2,e_3) \in E^3\mid e_1+e_2+e_3=0\}.$$
Clearly, $Y$ is an abelian surface over $\Q$ isomorphic to $E^2$.
Denote by $s$ an automorphism of $Y$ induced by the cyclic permutation of factors of $E^3$, i.e.,
$$s(e_1,e_2,e_3)=(e_3,e_1,e_2) \quad \forall\ (e_1,e_2,e_3) \in Y.$$
Let $C$ be the cyclic subgroup in $\Aut(X)$ of order $3$ generated by $s$.

By a theorem of Serre \cite{Serre2} the homomorphism
 $$\rho_{\ell,E}: G(\Q) \to \Aut_{\Z_{\ell}}(T_{\ell}(E)) \cong \GL(2,\Z_{\ell})$$ 
is {\sl surjective} for all but finitely many $\ell$.
Notice, that the composition
$$\det \rho_{\ell,E}: G(\Q) \to \GL(2,\Z_{\ell}) \to \Z_{\ell}^*$$
coincides with $\chi_{\ell}: G(K)\to \Z_{\ell}^*$ \cite{Serre2}.
In particular, if $\Q(E(\ell^{\infty}))$ is the field of definition of all
points on $E$ of $\ell$-power order
then $\Q(E(\ell^{\infty}))/\Q$ is the Galois extension with the Galois group
$\GL(2,\Z_{\ell})$. In addition, the cyclotomic
field $\Q(\ell)$ is the {\sl maximal abelian} subextension of $\Q(E(\ell^{\infty}))$ and the subgroup
$\Gal(\Q(E(\ell^{\infty}))/\Q(\ell)) \subset \Gal(\Q(E(\ell^{\infty}))/\Q)$ coincides with $\SL(2,\Z_{\ell})$.

Let us fix such an $\ell$, assuming in addition that $\ell-1$ is divisible
by $3$ but not divisible by $9$.
Let $\mu_{3,\ell}$ be the group of cubic roots of unity in $\Z_{\ell}^*$.
Then there exists a continuous
surjective homomorphism
$\pri_3: \Z_{\ell}^* \to \mu_{3,\ell},$
 coinciding with the identity map on $\mu_{3,\ell}$. 
These properties determine $\pri_3$ uniquely.

Let us define field $L$ as a subextension of $\Q(E(\ell^{\infty}))$ such that $\Q(E(\ell^{\infty}))/L$ is a cubic extension, whose Galois (sub)group coincides with
$$\mu_{3,\ell}\cdot\Id=\{\gamma \cdot \Id\mid \gamma \in \mu_{3,\ell}\} \subset \GL(2,\Z_{\ell}).$$
It follows immediately that $L$ is a Galois extension of $\Q$ and it does
not contain a primitive $\ell$th root of unity. 
 This implies  that $1$ and $-1$ are the only  roots of unity in $L$.

Let us choose a {\sl primitive} cubic root of unity $\gamma \in \mu_{3,\ell}$ and let
$\iota: \mu_{3,\ell} \to C$
be the isomorphism, which sends $\gamma$ to
$s$.

Now, let us define $X$ as the twist of $Y$ via the cubic character
$$\kappa:=\iota\pri_3 \chi_{\ell}
=\iota \pri_3 \det \rho_{\ell,E}
: G(\Q)\to \mu_{3,\ell} \to C \subset \Aut(Y).$$

The Galois module $T_{\ell}(X)$ is the twist of $T_{\ell}(E)^2$ via $\kappa$. Namely,
$$T_{\ell}(X)=\{(v_1,v_2,v_3)\in T_{\ell}(E)\oplus T_{\ell}(E)\oplus T_{\ell}(E)\mid v_1+v_2+v_3=0\}$$
as the $\Z_{\ell}-$module but
$$\rho_{\ell,X}(\sigma)(v_1,v_2,v_3)=\kappa(\sigma)(\rho_{\ell,E}(\sigma)(v_1),\rho_{\ell,X}(\sigma)(v_2),\rho_{\ell,X}(\sigma)(v_3))$$
for all $\sigma \in G(\Q)$. Now, we construct explicitly $G(L)-$invariant
elements of $T_{\ell}(X).$ Starting
with any $v \in T_{\ell}(E)$, put

$$w=(\gamma^{-1}v,\gamma v,v)=(\gamma^2 v,\gamma v, v)\in T_{\ell}(E)\oplus T_{\ell}(E)\oplus T_{\ell}(E).$$
Clearly, 
$w \in T_{\ell}(X);\quad sw =\gamma w.$
Let us check that $w$ is $G(L)-$invariant. Clearly,
$$G(L)=\{\sigma\in G(\Q)\mid \rho_{\ell,E}(\sigma) \in \mu_{3,\ell}\cdot \Id\}.$$
Let $\sigma \in G(L)$ with $\rho_{\ell,E}(\sigma)=\zeta\Id, \quad \zeta \in
\mu_{3,\ell}.$ If $\zeta=1$ , i.e., $\rho_{\ell,E}(\sigma)=\Id$ then all
elements of $V_{\ell}(X)$ are
$\sigma-$invariant. Since $\mu_{3,\ell}=\{1,\gamma,\gamma^{-1}\}$, we may
assume that $\zeta=\gamma$, i.e., 
$\rho_{\ell,E}(\sigma)=\gamma\cdot \Id$
 and therefore
$\det \rho_{\ell,E}(\sigma)=\gamma^2=\gamma^{-1}.$
Then
$$\rho_{\ell,X}(\sigma)(w)=$$
$$\iota(\pri_3(\det \rho_{\ell,E}(\sigma)))(\rho_{\ell,E}(\sigma)(\gamma^2
v),\rho_{\ell,E}(\sigma)(\gamma  v), \rho_{\ell,E}(\sigma)( v))=\iota(\gamma^2)(\gamma w)=$$
$$s^2(\gamma w)=\gamma s^2 w=\gamma\gamma^{2} w=w.$$

This proves that $w$ is $G(L)-$invariant.

Now, I claim that $X(L)$ contains infinitely many points, whose order is
a power of $\ell$. Indeed, starting with a non-divisible element $v \in T_{\ell}(E)$
and identifying the group $X_{\ell^n}$ with the quotient
$T_{\ell}(X)/\ell^n T_{\ell}(X)$, we get a $L-$rational point
$(\gamma^2 v,\gamma v,v)\mod\ell^n T_{\ell}(X)  \in T_{\ell}(X)/\ell^n
T_{\ell}(X)=X_{\ell^n}$ 
of order $\ell^n$.

\section{Another Example}

Let $K$ be an imaginary quadratic field with class number $1$ and let
$E$ be
an elliptic curve over $\Q$ such that $\End_K(E)=O_K$ is the ring of
integers in $K$. In this section we construct a Galois
extension $L$ of $K$ such that $E(L)$ contains 
infinitely many points of finite order but the intersection of $L$ and
$K(c)$  is of finite degree over $K$ (even coincides with $K$).

We write $\iota:\C \to \C$ for the complex conjugation $z\mapsto
\bar{z}$. We write $R$ for $O_K$. Clearly, $\End_{\Q}(E)=\Z\ne R$.
It follows easily that
$$\iota(ux)=\bar{u}(\iota(x))\quad \forall x \in E(\C), u\in R.$$
Notice that $K$ is abelian over $\Q$. Since $\Q(c)=\Q(ab)$, $K \subset
\Q(c)$ and therefore
$$K(c)=\Q(c).$$
\subsection{}
Let $\ell$ be a prime number. We write $R_{\ell}$ for $R \otimes \Z_{\ell}$.
It is well-known  that $T_{\ell}(E)$ is a free $R \otimes \Z_{\ell}$-module
of rank $1$ and therefore
$$\End_{R_{\ell}}(T_{\ell}(E))=R_{\ell},\quad \Aut_{R_{\ell}}(T_{\ell}(E))=R_{\ell}^*.$$
Let us consider the corresponding $\ell$-adic representation
$$\rho_{\ell,E}: G(\Q) \to \Aut_{\Z_{\ell}}(T_{\ell}(E)) \cong \GL(2,\Z_{\ell}).$$
Clearly, $G_{\ell}:=\rho_{\ell,E}(G(\Q))$ is not a subgroup of
$R_{\ell}^*=\Aut_{R_{\ell}}(T_{\ell}(E))$ but
$$H_{\ell}:=\rho_{\ell,E}(G(K))\subseteq R_{\ell}^*.$$

It is also known (\cite{Serre2}, Sect. 4.5) that
$$H_{\ell}=R_{\ell}^*$$
for all but finitely many primes $\ell$. Let us fix such an $\ell$, assuming
in addition that $\ell$ is unramified and splits in $K$.
This implies that $\ell=\q\bar{\q}$ for some  $\q\in K$ and
$$O_K=\q\cdot O_K+\bar{\q}\cdot O_K,\quad R_{\ell}=R_{\q}\oplus R_{\bar{\q}},
\quad R_{\q}=\Z_{\ell}, R_{\bar{\q}}= \Z_{\ell},$$
$$\q R_{\ell}=\ell R_{\q}\oplus R_{\bar{\q}}=
\ell \Z_{\ell}\oplus\Z_{\ell}\subset \Z_{\ell}\oplus\Z_{\ell}= R_{\ell},$$
$$\bar{\q} R_{\ell}= R_{\q}\oplus \ell R_{\bar{\q}}=
 \Z_{\ell}\oplus\ell\Z_{\ell}\subset \Z_{\ell}\oplus\Z_{\ell}= R_{\ell},$$
$$R_{\ell}^*=R_{\q}^*\times R_{\bar{\q}}^*,\quad R_{\q}^*=\Z_{\ell}^*, R_{\bar{\q}}^*= \Z_{\ell}^*.$$
We also have
$$T_{\ell}(E)=T_{\q}(E)\oplus T_{\bar{\q}}(E)$$
where
$$T_{\q}(E):=R_{\q}\cdot T_{\ell}(E),\quad T_{\bar{\q}}(E):=R_{\bar{\q}}\cdot T_{\ell}(E)$$
are free $\Z_{\ell}$-modules of rank $1$. This implies that for each positive integer $i$
$$\q^i T_{\q}(E)=\ell^i T_{\q}(E), \quad \bar{\q}^i T_{\q}(E)=T_{\q}(E),$$
$$\bar{q}^i T_{\bar{\q}}(E)=\ell^i T_{\bar{\q}}(E),\quad 
\q^i T_{\bar{\q}}(E)=T_{\bar{\q}}(E)$$
and therefore
$$T_{\ell}(E)/\ell^i T_{\ell}(E)=
T_{\q}(E)/\ell^i T_{\q}(E)\oplus T_{\bar{\q}}(E)/\ell^i T_{\bar{\q}}(E)=
T_{\q}(E)/\q^i T_{\q}(E)\oplus T_{\bar{\q}}(E)/\bar{\q}^i T_{\bar{\q}}(E).$$
It follows easily that a point
$x \in E_{\ell^i}=T_{\ell}(E)/\ell^i T_{\ell}(E)$ satisfies $\q^i x=0$ (respectively $\bar{\q}^i x=0$) if and only if
 $x \in T_{\q}(E)/\ell^i T_{\q}(E)$
(respectively $x \in T_{\bar{\q}}(E)/\ell^i T_{\bar{\q}}(E)$).

Let us put
$$\tau:=\rho_{\ell,E}(\iota) \in G_{\ell} \subset \Aut_{\Z_{\ell}}(T_{\ell}(E)).$$
Then $\tau^2=\Id$ and
$$\tau( R_{\q}^*\times\{1\}) \tau^{-1}=\{1\}\times R_{\bar{\q}}^*\subset R_{\ell}^* ,\quad
\tau(\{1\}\times R_{\bar{\q}}^*) \tau^{-1}=R_{\q}^*\times\{1\})\subset R_{\ell}^*.$$
It is also clear that
$$\tau(T_{\q}(E))=T_{\bar{\q}}(E), \quad \tau(T_{\bar{\q}}(E))=T_{\q}(E).$$

Let us consider the field $K(E(\ell^{\infty}))$  of definition of all
points on $E$ of $\ell$-power order. It is the Galois extension of $K$ with the Galois group
$R_{\ell}^*=R_{\q}^*\times R_{\bar{\q}}^*.$ 
It is also normal over $\Q$ and $\Gal(K(E(\ell^{\infty}))/\Q)=G_{\ell}$, since $E$ is defined over $\Q$ and $K$ is normal over $\Q$.

Let us define $L$ as a subextension of $K(E(\ell^{\infty}))/K$ such
that
$$\Gal(K(E(\ell^{\infty}))/L)=\{1\}\times R_{\bar{\q}}^*\subset
R_{\q}^*\times R_{\bar{\q}}^*=R_{\ell}^*=\Gal(K(E(\ell^{\infty}))/K).$$
One may
easily check that $L$ coincides with the field $K(E(\q^{\infty}))$  of definition of all
torsion points on $E$ which are killed by a power of $\q$. In
particular, $E(L)$ contains infinitely many points, whose order is a power
of $\ell$.
Let us consider the field $L'=\iota(L)$. Clearly, $K\subset L'\subset K(E(\ell^{\infty}))$
and $L'$ coincides with  the field $K(E(\bar{\q}^{\infty}))$
of definition of all torsion points on $E$ which are killed by a power of
$\bar{\q}$.
It is also clear that
$$\Gal(K(E(\ell^{\infty}))/L)=\tau(\{1\}\times
R_{\bar{\q}}^*)\tau^{-1}=R_{\q}^*\times\{1\}\subset
R_{\q}^*\times R_{\bar{\q}}^*=R_{\ell}^*=\Gal(K(E(\ell^{\infty}).$$

Since the subgroups $\{1\}\times R_{\bar{\q}}^*$ and
$R_{\q}^*\times\{1\}$ generate the whole group $R_{\ell}^*=\Gal(K(E(\ell^{\infty}))/K)$,
$$L\bigcap \iota(L)=L\bigcap L'=K.$$
It follows that if $M/K$ is a subextension of $L/K$ such that $M$ is
normal over $\Q$ then $M=K$. Since $K(c)=\Q(c)$,
$L\bigcap K(c)=L\bigcap \Q(c)$ is a subfield of
$\Q(c)$ and therefore is normal (even abelian) over $\Q$. It follows that
$$L\bigcap K(c)=K.$$

\section{Abelian subextensions}
The following statement may be viewed as a variant of Theorem \ref{Theorem 1} for
arbitrary abelian varieties over number fields.

\begin{thm}
\label{abelian}
\sl Let $X$ be an abelian variety over a number field $K$. Then:
\begin{enumerate}
\item  If for some prime $\ell$  the $\ell-$primary part of $\TORS(X(L))$
is infinite then $L$ contains an
abelian infinite subextension $E\subset L$ such that $\Gal(E/K)\cong \Z_{\ell}$ and $E/K$ is ramified only at divisors of $\ell$.

\item Let $P=P(X,L)$ be the set of primes $\ell$ such that  $X(L)$ contains
a point of order $\ell$.  If
$P$  is infinite then for all but finitely many primes $\ell \in P$ there
exist a finite  subextension
$E^{(\ell)}\subset L$ such that $E^{(\ell)}/K$ is a ramified abelian extension
which is unramified outside
divisors of $\ell$. 
In addition, the degree $[E^{(\ell)}:K]$ is prime to $\ell$ and degree
$[E^{(\ell)}:K]$ tends to infinity while $\ell$ tends to infinity.
\end{enumerate}
\end{thm}

\begin{cor}[Theorem of Bogomolov] If $\TORS(X(L))$  is infinite then $L$ contains an infinite abelian subextension of $K$.
\end{cor}
\begin{proof}[Proof of Theorem \ref{abelian}]
First, we may and will assume that $X$ is $K-$simple, i.e., the center
$F$ of the endomorphism algebra of $X$ is a number field.

Second, there is a positive integer $d$, enjoying the following property:

If $m$ is a positive integer such that $\varphi(m) \le 2g=2\dim(X)$ then
$d$ is divisible by $m$.

Third, let $\lambda$ be a prime ideal in $O_F$ dividing a prime number
$\ell$. Then, in the notations of Section 1 the following statement is
true.
\begin{lem}
\begin{enumerate}
\item
The composition
$$\pi_{\lambda}:=(\det_{F_{\lambda}} \rho_{\lambda,X})^d:G(K) \to
\Aut_{F_{\lambda}}V_{\lambda}(X) \to {F_{\lambda}}^* \to
{F_{\lambda}}^*$$
is an abelian representation of $G(K)$ unramified outside divisors of
$\ell$.

\item For all but finitely many $\lambda$ the composition
$$\bar{\pi}_{\lambda}:=(\det_{F_{\lambda}}\bar{\rho}_{\lambda,X})^d:G(K) \to
\Aut_{O_F/\lambda}X_{\lambda} \to (O_F/\lambda)^* \to
(O_F/\lambda)^*$$
is an abelian representation of $G(K)$ unramified outside divisors of
$\ell$.

\end{enumerate}
\end{lem}

We will prove Lemma at the end of this section. Now, let us finish the
proof of Theorem, assuming validity of Lemma.

First, notice that the ratio 
$$e=2\dim(X)/[F:\Q]$$
is a positive integer. Second,
for all but finitely many primes $p$  there exists a finite
collection  of {\sl Weil numbers}, i.e., certain algebraic integers  $\{\alpha_1,\ldots
\alpha_{e}\} \subset F(a)$, enjoying the following properties:
\begin{itemize}
\item (Weil's condition) There is a positive integer $q>1$ such that $q$ is an integral power
of $p$ and all $\mid\alpha_i\mid^2=q$ for all embeddings
$F(a)\subset\C$. 
\item
For all $\ell\ne p$ and $\lambda\mid\ell$ there is a subset
$S_{\lambda}\subset \{1,\ldots e\}$ such that
$(\prod_{i\in S_{\lambda}}\alpha_i)\in
O_F$ and the group $\Im(\pi_{\lambda})$ contains $\prod_{i\in
S_{\lambda}}\alpha_i$.
\item For all but finitely many $\lambda$ the subgroup
$\Im(\bar{\pi}_{\lambda})$
contains $(\prod_{i\in S_{\lambda}}\alpha_i)\mod \lambda \in (O_F/\lambda)^*$.
\end{itemize}

Indeed, let us choose  a prime ideal $\mathbf v$ in the ring $O_K$ of all
algebraic integers in $K$ such that $X$ has good reduction at $\mathbf
v$. Let
$$Fr_{\mathbf v} \in \Im(\rho_{\lambda,X}) \subset \Aut_{F_{\lambda}}V_{\lambda}(X)$$
be {\sl Frobenius element} $Fr_{\mathbf v}$ at $\mathbf v$ (defined up to
conjugacy)\cite{Serre},\cite{RibetA}. Then the set of its eigenvalues belongs
to $F(a)$, does not depend on the choice of $\lambda$ and satisfies all the desired properties with $p$
the residual characteristic of $\mathbf v$ and $q=\#(O_K/{\mathbf v})$(\cite{Shimura},
Ch. 7, Prop. 7.21 and proof of Prop. 7.23).

\begin{proof}[Proof of assertion 1] We know that there exists $\lambda$ dividing
$\ell$ such that the subspace $V_{\lambda}(X)$ consists of $G(L)-$invariants. This means
that $G(L)$ lies in the kernel of $\pi_{\lambda}$. This implies that
the field $E'$ of $\ker(\pi_{\lambda})-$invariants is an abelian subextension
of L, unramified outside divisors of $\ell$ and $\Gal(E'/K)$ is
isomorphic to $\Im(\pi_{\lambda})$. Choosing a collection of Weil numbers
attached to prime $p\ne\ell$, we easily conclude that $\Im(\pi_{\lambda})$ is an
{\sl infinite} commutative $\ell-$adic Lie group \cite{Serre} and therefore,
there is a continuous quotient of $\Im(\pi_{\lambda})$, isomorphic to
$\Z_{\ell}$. One has to take as $E$ the subextension of $E'$
corresponding to this quotient.
\end{proof}
\begin{proof}[Proof of assertion 2]
We know that for all but finitely many $\ell \in P$ there exists
$\lambda$ dividing $\ell$ such that $X_{\lambda}$ consists of $G(L)-$invariants.
This means that the field $E^{(\ell)}$ of 
$\ker(\bar{\pi}_{\lambda})-$invariants is an abelian subextension
of L, unramified outside divisors of $\ell$ and $\Gal(E^{(\ell)}/K)$ is
isomorphic to $\Im(\pi_{\lambda})$. In order to prove that
$[E^{(\ell)}:K]$ tends to infinity, let us assume that there exist an
infinite subset $P'\subset P$ and a positive integer $D$ such that
$\#(\Gal(E^{(\ell)}/K))=[E^{(\ell)}:K]$ divides $D$ for all $\ell\in P'$. This means that
$$\bar{\pi}_{\lambda}^D: G(K)\to (O_F/\lambda)^*$$
is a trivial homomorphism for {\sl infinitely many} $\lambda$. 
In order to get a contradiction, let us choose a collection of Weil
numbers $\{\alpha_1,\ldots
\alpha_{e}\}$ enjoying the properties described above.

Clearly. for any non-empty subset $ S \subset \{1,\ldots e\}$ the
product $\alpha_S:=\prod_{i\in S} \alpha_i$ is not a root of unity. In addition,
if $\alpha_S\in
O_F$ then there only finitely many $\lambda$ such that $\alpha_S^D-1$ is
an element of $\lambda$. Since there are only finitely many subsets of $\{1,\ldots2g\}$,
for all but finitely many $\lambda$ the group
$$\Im((\bar{\pi}_{\lambda})^D) \subset (O_F/\lambda)^*$$
contains an element of type $\alpha_S^D\mod \lambda$ different from 1.
This implies that $(\bar{\pi}_{\lambda})^D$ is a non-trivial
homomorphism for all but finitely many $\lambda$. This gives the
desired contradiction.

\end{proof}
\begin{proof}[Proof of Lemma]
Let  $\mathbf v$ be a prime ideal in the ring $O_K$ of all
algebraic integers in $K$. We write $I_{\mathbf v} \subset G(K)$ for
the corresponding inertia subgroup defined up to conjugacy. Assume that the
residual characteristic of $\mathbf v$ is different from $\ell$. It is known
\cite{SGA} that for any $ \sigma \in I_{\mathbf v}$ there exists a
positive integer $m$ such that $\rho_{\ell,X}(\sigma)^m$ is an unipotent
operator in $V_{\ell}(X)$ and its characteristic polynomial has
coefficients in $\Z$. This implies that if $m$ is the smallest integer enjoying
this property then the characteristic polynomial is divisible by the
$m$th cyclotomic polynomial. This implies that $2g\ge \varphi(m)$ and
therefore $m$ divides $d$. Since $V_{\lambda}(X)$ is a Galois-invariant
subspace of $V_{\ell}(X)$ and (for all but finitely many $\ell$) $X_{\lambda}$
is a Galois-invariant
subspace of $T_{\ell}(X)/\ell T_{\ell}(X)$, a Galois automorphism
$\sigma^d$ acts as an unipotent operator in $V_{\lambda}(X)$ and (for
all but finitely many $\lambda$) in $X_{\lambda}$. One has only to
recall that the determinant of an unipotent operator is always $1$.
\end{proof}
\end{proof}

\end{document}